\documentclass[12pt]{amsart}

\usepackage{latexsym,amssymb} 
\usepackage{amsmath,amsfonts,amsthm}
\input xypic

\newtheorem{theorem}{Theorem}[section]
\newtheorem{proposition}[theorem]{Proposition}
\newtheorem{corollary}[theorem]{Corollary}
\newtheorem{lemma}[theorem]{Lemma}
\newtheorem{example}[theorem]{Example}

\begin{document}

\title[Compactable semilattices]{Compactable semilattices}

\author{M.~V.~Lawson}

 \address{Department of Mathematics and the
Maxwell Institute for Mathematical Sciences\\
Heriot-Watt University\\
Riccarton\\
Edinburgh~EH14~4AS\\
Scotland}

\email{M.V.Lawson@ma.hw.ac.uk}

\thanks{This research was supported by an EPSRC grant (EP/F004184, EP/F014945, EP/F005881).
I am grateful to Daniel Lenz and Pedro Resende for email correspondence on this paper.}

\dedicatory{This paper is dedicated to the memory of my colleague Prof.~W.~Douglas~Munn}



\begin{abstract} 
We characterize those semilattices that give rise to Boolean spaces on their associated spaces of ultrafilters.
The class of $0$-disjunctive semilattices, important in the theory of congruence-free inverse semigroups,
plays a distinguished role in this theory.
\end{abstract}

\maketitle

\section{Introduction}

In Stone duality, a topological space is associated with each Boolean algebra $B$, the elements of this space being the ultrafilters on $B$.
Such a space is known as a Boolean space and is compact and hausdorff with a base of clopen sets.
Furthermore, the set of clopen subsets of this space is then a Boolean algebra isomorphic to $B$ \cite{BS}.

The goal of this paper is to characterize those bounded semilattices whose spaces of ultrafilters are Boolean.
This question arose from work generalizing \cite{Lawson}
and from the diagram at the bottom of page 471 of \cite{M1},
and seemed to us to be of independent interest.

The solutions to finding our characterizations make essential use of ideas due to Exel \cite{Exel}, in our Theorem~2.5,
and Lenz \cite{Lenz}, in our Theorem~2.10.
For basic topology see \cite{S,W}.

\section{The main theorems}

Throughout this paper we will deal only with bounded meet semilattices $E$:
that is, semilattices with both a top and a bottom.
Elements $e,f \in E$ are said to be {\em orthogonal} if $e \wedge f = 0$.
We denote by $e^{\ast}$ the set of all elements of $E$ orthogonal to $e$.

If $X \subseteq E$, define 
$$X^{\uparrow} = \{e \in E \colon f \leq e \text{ for some } f \in X \}$$
and
$$X^{\downarrow} = \{e \in E \colon e \leq f \text{ for some } f \in X \}.$$
 A subset $X \subseteq E$ is called an {\em order ideal} if $X = X^{\downarrow}$.
The {\em principal order ideal generated by $e$} is $e^{\downarrow}$.

A subset $F$ of $E$ is called a {\em filter} if it satisfies the following three conditions:
\begin{description} 

\item[{\rm (F1)}] $0 \notin F$.

\item[{\rm (F2)}] $e,f \in F$ implies that $e \wedge f \in F$.

\item[{\rm (F3)}] $e \in F$ and $e \leq f$ implies that $f \in F$.
 
\end{description}
If $F$ is a filter then $F^{\uparrow} = F$.

An {\em ultrafilter} is a maximal filter.

Each non-zero element of $E$ is contained in a filter: namely, $e^{\uparrow}$, the {\em principal filter} containing $e$.
Thus the set of all filters containing $e$ is non-empty.
This set is ordered by inclusion and has the property that the union of every chain is again a filter containing $e$.
The following is therefore a consequence of Zorn's Lemma.

\begin{lemma} 
Every non-zero element of a semilattice is contained in an ultrafilter.
\end{lemma}

The following very useful result is Lemma~3.2 of \cite{Exel}.

\begin{lemma} Let $F$ be a filter in a semilattice $E$.
Then $F$ is an ultrafilter if and only if $F$ contains every element $b \in S$ such that $b \wedge c \neq 0$ for all $c \in F$.
\end{lemma}

Denote by $\mathsf{U}(E)$ the set of all ultrafilters on $E$.
For each $e \in E$ define $\mathcal{K}_{e}$ to be the set of all ultrafilters that contain $e$.
Clearly $\mathcal{K}_{0} = \emptyset$, because no ultrafilter contains 0,
and this is the only empty such set,
and $\mathcal{K}_{1} = \mathsf{U}(E)$, because every ultrafilter contains 1.
Put 
$$\Omega = \Omega (E) = \{\mathcal{K}_{e} \colon e \in E \}.$$

\begin{lemma} Let $E$ be a semilattice.
Then $\Omega$ is a base of clopen sets for a hausdorff topology on the set $\mathsf{U}(E)$. 
\end{lemma}
\begin{proof}
It is easy to check that $\mathcal{K}_{e} \cap \mathcal{K}_{f} = \mathcal{K}_{e \wedge f}$ for all $e,f \in E$
and so $\Omega$ is a base for a topology on $E$.

We now prove that $\mathsf{U}(E) \setminus \mathcal{K}_{e}$ is an open set which implies that
each of the sets $\mathcal{K}_{e}$ is also closed.
Let $F \in \mathsf{U}(E) \setminus \mathcal{K}_{e}$.
By assumption $e \notin F$.
Thus by Lemma~2.2, there exists $f \in F$ such that $e \wedge f = 0$.
It follows that $F \in \mathcal{K}_{f}$.
But $F \in \mathcal{K}_{f}$ implies that $e \notin F$.
Thus $F \in \mathcal{K}_{f} \subseteq \mathsf{U}(E) \setminus \mathcal{K}_{e}$.
It follows that $\mathsf{U}(E) \setminus \mathcal{K}_{e}$ is open and so $\mathcal{K}_{e}$ is closed.
We have therefore shown that $\mathsf{U}(E)$ has a base of clopen sets.

Let $F,G \in \mathsf{U}(E)$ where $F \neq G$.
There exists $e \in F$ such that $e \notin G$.
It follows by Lemma~2.2 that there exists $f \in G$ such that $e \wedge f = 0$.
Hence $\mathcal{K}_{e} \cap \mathcal{K}_{f} = \emptyset$
and $F \in \mathcal{K}_{e}$ and $G \in \mathcal{K}_{f}$.
We have therefore shown that $\mathsf{U}(E)$ is hausdorff.
\end{proof}

In this paper, we shall always regard $\mathsf{U}(E)$ as a topological space with respect to the specific base $\Omega(E)$.
We denote the topology it generates by $\tau$.
We say that a bounded meet semilattice $E$ is {\em compactable} if the topological space $(\mathsf{U}(E),\tau)$ is compact.
It follows that the spaces associated with compactable semilattices are boolean.

We shall now use the ideas of Exel \cite{Exel}.
Let $X,Y \subseteq E$.
The subset $E^{X,Y}$ of $E$ consists of all those elements of $E$ that are below
every element of $X$ and orthogonal to every element of $Y$.
We allow the possibility that $Y$ might be empty.\\

\noindent
{\bf Remark }In practice, $X$ will consist of one element and $Y$ will be finite.
Thus we shall be interested in sets of the form
$e^{\downarrow} \cap \{e_{1}, \ldots, e_{m} \}^{\ast}$.\\

A subset $Z \subseteq E^{X,Y}$ is said to be a {\em cover}
if for each non-zero $e \in E^{X,Y}$ there exists $z \in Z$ such that $e \wedge z \neq 0$.
Thus a cover is a `large' subset in the sense that every element of  $E^{X,Y}$
meets an element of $Z$.

A filter $F \subseteq E$ is said to be {\em tight} if the following condition is satisfied:
for all finite sets $X$ and $Y$
\begin{center}
$X \subseteq F$
and
$Y \cap F = \emptyset$
and $Z$ a finite cover of $E^{X,Y}$ $\Rightarrow$  $Z \cap F \neq \emptyset$.
\end{center}

The following is Proposition~4.2 of \cite{Exel}.
We give the proof for the sake of completeness.

\begin{lemma} Every ultrafilter is tight.
\end{lemma}
\begin{proof} Let $F$ be an ultrafilter in $E$.
Let $X$ and $Y$ be finite sets such that $X \subseteq F$ and $Y \cap F = \emptyset$.
Let $Z$ be a finite cover of $E^{X,Y}$.
We shall prove that $Z \cap F \neq \emptyset$.
Suppose not.
Then $Z \cap F = \emptyset$.
By Lemma~2.2, for each $z \in Z$ there exists $x_{z} \in F$ such that $z \wedge x_{z} = 0$,
and for each $y \in Y$ there exists $x_{y} \in F$ such that $y \wedge x_{y} = 0$.
Put
$$ w = \left( \bigwedge_{x \in X} x \right) \wedge \left( \bigwedge_{y \in Y} x_{y}  \right) \wedge \left( \bigwedge_{z \in Z} x_{z} \right).$$ 
Then $w \in F$ by construction and so $w \neq 0$.
But clearly $w \in E^{X,Y}$.
Thus there exists $z_{1} \in Z$ such that $w \wedge z_{1} \neq 0$.
But $w \leq x_{z_{1}}$ and $x_{z_{1}} \wedge z_{1} = 0$ implies that $w \wedge z_{1} = 0$.  
This is a contradiction.
Observe that the above argument also applies to the case where $Y$ is empty.\end{proof}

The following is really an easy deduction from \cite{Exel} and \cite{Lenz} but we give a complete proof because of its importance:
it is our first characterization of compactable semilattices.

\begin{theorem} A semilattice is compactable if and only if every tight filter is an ultrafilter.
\end{theorem}
\begin{proof} 
Suppose that $E$ is compactable.
By Lemma~2.4, it is enough to prove that every tight filter is an ultrafilter.
Let $F$ be a tight filter.
Let $e \in E$ be an element such that $e \wedge f \neq 0$ for all $f \in F$.
We shall show that $e \in F$ which will prove that $F$ is an ultrafilter by Lemma~2.2.
There are two cases to consider: $e^{\ast} = 0$ and $e^{\ast} \neq 0$.

Suppose first that $e^{\ast} = 0$.
Let $f \in F$ be arbitrary.
Put $X = \{f\}$ and $Y = \emptyset$.
Then $X \subseteq F$. 
Put $e' = e \wedge f \neq 0$.
Then $e' \leq f$.
It follows that $e' \in E^{X,Y}$.
Let $j \in E^{X,Y}$ be an arbitrary non-zero element.
Then $j \leq f$.
Observe that $j \wedge e' = j \wedge e \wedge f = j \wedge e \neq 0$ by assumption that $e^{\ast} = 0$.
It follows that $e'$ is a cover for $E^{X,Y}$.
By our assumption that $F$ is tight it follows that $e' \in F$ and so $e \in F$, as required.

Suppose now that $e^{\ast} \neq 0$.
We claim that there exists a finite subset $\{e_{1}, \ldots, e_{n} \} \subseteq e^{\ast}$
such that  
$$\mathsf{G}(E) \setminus \mathcal{K}_{e} =  \bigcup_{i=1}^{n} \mathcal{K}_{e_{i}}.$$
To see why, let $G$ be an arbitrary ultrafilter such that $e \notin G$. 
Then there exists $g \in G$ such that $e \wedge g = 0$.
It follows that
$$\mathsf{G}(E)\setminus \mathcal{K}_{e} = \bigcup_{i \in e^{\ast}} \mathcal{K}_{i}.$$
This is an open cover of a closed subset of a compact space and so an open cover of a compact subset.
It follows that there is a finite subcover.
There are therefore a finite number of non-zero elements $e_{1}, \ldots, e_{n}$ which are orthogonal to $e$,
that is $Y = \{e_{1}, \ldots, e_{n}\} \subseteq e^{\ast}$,
such that
$$\mathsf{G}(E) \setminus \mathcal{K}_{e} = \bigcup_{i=1}^{n} \mathcal{K}_{e_{i}}$$
and the claim is verified.

Now let $f \in F$ be arbitrary.
Put $X = \{f\}$ and $Y = \{e_{1}, \ldots ,e_{n} \}$, where $Y$ is the set chosen above.
Then $X \subseteq F$. 
Also $Y \cap F = \emptyset$ because if $e_{i} \in F$ then $e \wedge e_{i} = 0$ contradicts our assumption about $e$.
Put $e' = e \wedge f \neq 0$.
Then $e' \leq f$ and $e' \wedge e_{i} = e \wedge f \wedge e_{i} = 0$.
It follows that $e' \in E^{X,Y}$.
Let $j \in E^{X,Y}$ be an arbitrary non-zero element.
Then $j \leq f$ and $j \wedge e_{i} = 0$ for all $i$.
Now $j \wedge e' = j \wedge e \wedge f = j \wedge e$.
Let $G$ be any ultrafilter containing $j$.
Then $e \in G$ because if not $G$  must contain at least one of the $e_{i}$.
But then $j,e_{i} \in G$ implies that $j \wedge e_{i} \neq 0$ contradicting our choice of $j$.
Thus $j \wedge e \neq 0$.
It follows that $e'$ is a cover for $E^{X,Y}$.
By our assumption that $F$ is tight it follows that $e' \in F$ and so $e \in F$, as required.

Suppose now that every tight filter is an ultrafilter.
We prove that $E$ is compactable which means that we have to prove that $(\mathsf{U}(E),\tau)$ is compact.
This will be proved in stages.

We denote by $\mathbf{2}^{E}$ the set of all functions from $E$ to the set $\mathbf{2} = \{0,1 \}$,
regarded as the semilattice $0 < 1$.
This is the product of the set $\mathbf{2}$ with itself $\left|E \right|$ times.
Giving $\mathbf{2}$ the discrete topology, we see by Tychonoff's theorem that  $\mathbf{2}^{E}$ is compact.
A subbase for this topology is given by subsets of the form
$U_{e}$ and $U_{e}^{c}$ where $e \in E$ and 
$$U_{e} = \{\theta \colon E \rightarrow \mathbf{2} \colon \theta(e) = 1 \}
\text{ and }
U_{e}^{c} = \{\theta \colon E \rightarrow \mathbf{2} \colon \theta(e) = 0 \}.$$
These are clopen sets.

A function $\theta \colon E \rightarrow \mathbf{2}$ is called a {\em homomorphism} iff $\theta (e \wedge f) = \theta (e)\theta(f)$
for all $e,f \in E$.
Thus a function $\phi$ fails to be a homomorphism iff $\phi (e \wedge f) \neq \phi (e)\phi (f)$ for some $e,f \in E$
iff $\phi$ belongs to the union of sets of the form
$$U_{e} \cap U_{f} \cap U_{e \wedge f}^{c}
\text{ or }
U_{e}^{c} \cap U_{f} \cap U_{e \wedge f}
\text{ or }
U_{e} \cap U_{f}^{c} \cap U_{e \wedge f}
\text{ or }
U_{e}^{c} \cap U_{f}^{c} \cap U_{e \wedge f}$$
for some $e,f \in E$.
It follows that the set of functions which are not homomorphisms is open
and so the set of homomorphisms is closed.
A homomorphism $\theta \colon E \rightarrow \mathbf{2}$ is called a {\em representation} if $\theta(0) = 0$ and $\theta(1) = 1$. 
It follows that a function is a representation precisely when it lies in the intersection of three closed sets.
Thus the set of representations is a closed subset of a compact space and so forms a compact space.

There is a bijection between representations of $E$ and filters in $E$ \cite{P}:
if $F$ is a filter in $E$ define $\theta_{F}(e)$ to be the characteristic function of $F$;
conversely, if $\theta \colon E \rightarrow \mathbf{2}$ is a representation define $F = \theta^{-1}(1)$, which is a filter.
The topology induced by this bijection on the set $\mathsf{F}(E)$ of all filters on $E$
has as a base the sets of the form
$\mathcal{U}_{e:e_{1}, \ldots, e_{m}}$ where $e_{1}, \ldots, e_{m} \leq e$
and $\mathcal{U}_{e:e_{1}, \ldots, e_{m}}$ consists of all filters that contain $e$ and omit the $e_{i}$.
Thus $\mathsf{F}(E)$ is a compact space.
Now $\mathsf{U}(E)$ is a subset of  $\mathsf{F}(E)$.
We prove that the subspace topology agrees with $\tau$.\footnote{This could be deduced immediately from Lemma~6.8 of \cite{Lenz}.} 
Let $F$ be an ultrafilter in  $\mathcal{U}_{e:e_{1}, \ldots, e_{m}}$.
Then $e \in F$ and $\{e_{1}, \ldots, e_{m} \} \cap F = \emptyset$.
By Lemma~2.2, for each $i$ there exists $f_{i} \in F$ such that $e_{i} \wedge f_{i} = 0$.
Put $i = e \wedge f_{1} \wedge \ldots \wedge f_{m}$.
Then $i \in F$.
Consider $\mathcal{K}_{i}$.
Any ultrafilter in  $\mathcal{K}_{i}$ contains $e$ and omits all the $e_{i}$.
Thus $F \in \mathcal{K}_{i} \subseteq \mathcal{U}_{e:e_{1}, \ldots, e_{m}}$. 
Exel proves in Theorem~4.4 of \cite{Exel} that the closure of $\mathsf{U}(E)$ in $\mathsf{F}(E)$ is precisely the set of tight filters.
It follows under our assumption that $\mathsf{U}(E)$ is a closed subset of $\mathsf{F}(E)$ and so is compact.\end{proof}

Our goal now is to obtain a criterion for compactability that uses only ultrafilters.
This will necessarily involve restricting the class of semilattices we consider.

Let $f \neq 0$.
Following Lenz \cite{Lenz}, 
write 
$$f \rightarrow (e_{1}, \ldots, e_{m})$$
iff
for each $0 \neq x \leq f$ there exists $i$ such that $x \wedge e_{i} \neq 0$.
The following is Proposition~6.2 of \cite{Lenz}.
We give the proof for the sake of completeness.

\begin{lemma} Let $E$ be a semilattice.
Then $f \rightarrow (e_{1}, \ldots, e_{m})$
if and only if 
$\mathcal{K}_{f} \subseteq \bigcup_{i=1}^{m} \mathcal{K}_{e_{i}}$.
\end{lemma}
\begin{proof}
Suppose that $\mathcal{K}_{f} \subseteq \bigcup_{i=1}^{m} \mathcal{K}_{e_{i}}$.
Let $0 \neq i \leq f$. 
By Lemma~2.1, there is an ultrafilter $F$ containing $i$.
Thus $F$ contains $f$.
But by assumption $e_{i} \in F$ for some $i$.
It follows that $i \wedge e_{i} \neq 0$.
We have proved that $f \rightarrow (e_{1}, \ldots, e_{m})$.

To prove the converse, suppose that $f \rightarrow (e_{1}, \ldots, e_{m})$.
Let $F$ be an ultrafilter that contains $f$.
Suppose that $e_{i} \notin F$ for all $i$.
By Lemma~2.2, for each $i$ there exists $f_{i} \in F$ such that $e_{i} \wedge f_{i} = 0$.
Put $j = f_{1} \wedge f_{2} \wedge \ldots \wedge f_{m} \in F$.
Then $j \wedge e_{i} = 0$ for all $i$.
Now $j,f \in F$ and so $j \wedge f \neq 0$.
But from $0 \neq j \wedge f \leq f$ we deduce that $j \wedge f \wedge e_{i} \neq 0$ for some $i$
since $f \rightarrow (e_{1}, \ldots, e_{m})$.
Thus $j \wedge e_{i} \neq 0$ for some $i$, which is a contradiction.
Thus $e_{i} \in F$ for some $i$, as required.
\end{proof}

The above lemma is important because it provides a link between properties of elements and properties of ultrafilters.

\begin{corollary} Let $e$ and $f$ be non-zero elements of the semilattice $E$.
Then $\mathcal{K}_{e} \subseteq \mathcal{K}_{f}$ if and only if
$e \rightarrow f$.
\end{corollary}

A semilattice $E$ is said to be {\em $0$-disjunctive} if for all $e,f \in E\setminus \{0 \}$
such that $e < f$, there exists $0 \neq e' \leq f$ such that $e \wedge e' = 0$;
that is, $e^{\downarrow} \cap f^{\ast} \neq 0$.
Such semilattices implicitly arise in Exel's paper \cite{Exel} in the negative sense via the concept of one idempotent being dense in another
but they have a wider significance: the semilattices of idempotents of congruence-free
inverse semigroups are $0$-disjunctive \cite{M1, Petrich}.
Munn \cite{M1} points out that a semilattice with zero $E$ is $0$-disjunctive iff the following condition is satisfied:
for all $e,f \in E$ if $e \neq f$ then there exists $g \in E$ such that 
either $e \wedge g = 0, f \wedge g \neq 0$ or $e \wedge g \neq 0, f \wedge g \neq 0$.

We say that a semilattice is {\em separative} if and only if $\mathcal{K}_{e} = \mathcal{K}_{f}$ implies $e = f$.

\begin{lemma} Let $E$ be a semilattice.
Then the following are equivalent.
\begin{enumerate}

\item $E$ is $0$-disjunctive.

\item For all $e,f \in E \setminus \{ 0\}$ we have that $e < f \Rightarrow \mathcal{K}_{e} \subset \mathcal{K}_{f}$. 

\item $E$ is separative.

\end{enumerate}
\end{lemma}
\begin{proof} (1)$\Rightarrow$(2).
Suppose that $E$ is $0$-disjunctive.
Let $0 \neq e < f$.
This immediately implies that $\mathcal{K}_{e} \subset \mathcal{K}_{f}$
so our problem is to show that this is a strict inclusion.
By assumption there exists $0 \neq e'$ such that $e' \leq f$ and $e \wedge e' = 0$.
By Lemma~2.1, let $F$ be an ultrafilter containing $e'$.
Then $f \in F$ but clearly $e \notin F$.
It follows that $F \in \mathcal{K}_{f} \setminus \mathcal{K}_{e}$ and so $\mathcal{K}_{e} \subset \mathcal{K}_{f}$.

(2)$\Rightarrow$(3). Let $\mathcal{K}_{e} = \mathcal{K}_{f}$.
Now 
$\mathcal{K}_{e} \cap \mathcal{K}_{f} = \mathcal{K}_{e \wedge f}$.
Thus in this particular case we have that 
$\mathcal{K}_{e} 
= \mathcal{K}_{e \wedge f}
=
\mathcal{K}_{f}$.
If $e \wedge f < e$ then $\mathcal{K}_{e \wedge f} \subset \mathcal{K}_{e}$.
But we have equality of sets.
Thus $e \wedge f = e$.
Similarly $e \wedge f = f$.
It follows that $e = f$, as required.

(3)$\Rightarrow$(1). Let $0 \neq e < f$.
Then by assumption there is an ultrafilter $F$ containing $f$ and omitting $e$.
Thus by Lemma~2.2 there exists $f' \in F$ such that $f' \wedge e = 0$.
Put $i = f \wedge f' \in F$.
Then $0 \neq i \leq f$ and $i \wedge e = 0$.
Thus $E$ is $0$-disjunctive.
\end{proof}

\begin{lemma} Let $E$ be a compactable semilattice.
Suppose that $0 \neq f < e$ and that $e^{\downarrow} \cap f^{\ast} \neq 0$.
Then we can find finitely many elements $e_{1}, \ldots, e_{m}$ such that
$$e \rightarrow (e_{1}, \ldots, e_{m},f)
\text{ and }
e_{i} \in e^{\downarrow} \cap f^{\ast}.$$
\end{lemma}
\begin{proof} By assumption $\mathcal{K}_{f} \subseteq \mathcal{K}_{e}$ and $\mathcal{K}_{e} \setminus \mathcal{K}_{f} \neq \emptyset$. 
We may write
$$\mathcal{K}_{e} \setminus \mathcal{K}_{f}
=
\bigcup_{i \in e^{\downarrow} \cap f^{\ast}} \mathcal{K}_{i}.$$
To see why, let $F$ be an ultrafilter that contains $e$ and omits $f$.
By Lemma~2.2, there exists $g \in F$ such that $f \wedge g = 0$.
But $e,g \in F$ implies that $e \wedge g \in F$.
Thus $e \wedge g \leq e$ and $e \wedge g$ is orthogonal to $f$
and so $F$ also belongs to the righthand side.
It is clear that the reverse inclusion holds.
The set $\mathcal{K}_{e} \setminus \mathcal{K}_{f}$ is the intersection of two closed sets and so is closed.
Since $E$ is compactable it follows that $\mathcal{K}_{e} \setminus \mathcal{K}_{f}$ is compact.
We may therefore find a finite number of elements $e_{1}, \ldots, e_{m} \leq e$ such that
each $e_{i}$ is orthogonal to $e$ and 
$$\mathcal{K}_{e} \setminus \mathcal{K}_{f}
=
\bigcup_{i=1}^{m} \mathcal{K}_{e_{i}}.$$
Clearly, 
$$\mathcal{K}_{e} = \mathcal{K}_{f} \cup \bigcup_{e_i \leq e, e_i \in f^{\ast}} \mathcal{K}_{e_i}.$$
Thus by Lemma~2.6, we have that $e \rightarrow (e_{1}, \ldots, e_{m},f)$.
\end{proof}

A semilattice $E$ is said to satisfy the {\em trapping condition}
if for all $e,f \in E \setminus \{ 0\}$ where $f < e$ there exists $e_{1}, \ldots, e_{m}$, where $m \geq 1$, such that
$e \rightarrow (e_{1}, \ldots, e_{m},f)$
and
$e_{i} \in e^{\downarrow} \cap f^{\ast}$.\\

\noindent
{\bf Remarks }\mbox{}
\begin{enumerate}

\item We have phrased the trapping condition in such a way that if it holds the semilattice is $0$-disjunctive.

\item A trapping condition was introduced by Lenz \cite{Lenz} and ours is a special case of his.

\end{enumerate}

We may now give a characterization of an important class of compactable semilattices
which uses only ultrafilters and not arbitrary filters.

\begin{theorem} A separative semilattice is compactable if and only if it satisfies the trapping condition. 
\end{theorem}
\begin{proof} It follows from Lemma~2.9 that every separative compactable semilattice satisfies the trapping condition.
We prove the converse.

Let $E$ be a semilattice that satisfies the trapping condition.
By Theorem~2.5, we need to prove that every tight filter is an ultrafilter.
Let $F$ be a tight filter and assume that it is not an ultrafilter. 
Then we can find an ultrafilter $G$ such that $F \subset G$.
Let $g \in G \setminus F$ and let $f \in F$ be arbitrary.
Then $g' = g \wedge f \in G$.
We shall prove that the tightness condition on $F$ implies that $g' \in F$ giving $g \in F$ and so a contradiction.
We have that $0 \neq g' < f$.
Thus by the trapping condition there exist idempotents $f_{1}, \ldots, f_{n} \leq f$ 
such that $f \rightarrow (f_{1}, \ldots, f_{n},g')$ and $f_{i} \wedge g' = 0$.
Suppose that $f_{i} \in F$ for some $i$.
Then  $f_{i},g' \in G$ and so $f_{i} \wedge g' \neq 0$
but this contradicts the fact that $f_{i} \wedge g' = 0$. 
Thus $\{f_{1}, \ldots, f_{n} \} \cap F = \emptyset$.
We consider the set $E^{\{ f\}, \{f_{1}, \ldots, f_{n} \}}$ with respect to the tight filter $F$.
Clearly $g' \in E^{\{f\}, \{f_{1}, \ldots, f_{n} \} }$.
Let $0 \neq i \in  E^{\{f\}, \{f_{1}, \ldots, f_{n} \} }$.
From $0 \neq i \leq f$ and  $f \rightarrow (f_{1}, \ldots, f_{n},g')$ there are two possibilties.
Either $i \wedge f_{i} \neq 0$ for some $i$, which cannot happen by our choice of $i$,
or $i \wedge g' \neq 0$.
It follows that $g'$ is a cover for $E^{\{f\}, \{f_{1}, \ldots, f_{n} \} }$.
But $F$ is a tight filter and so $g' \in F$ which is a contradiction.
It follows that $F$ is an ultrafilter.
\end{proof}

The proof that the trapping condition implies compactability was proved in Proposition~6.7 of \cite{Lenz}.

If $E$ is a compactable semilattice then there is a boolean algebra $\mathfrak{B}(E) = \mathsf{B}(\mathsf{U}(E))$ 
whose elements are the clopen subsets of $\mathsf{U}(E)$.
The function $\kappa \colon E \rightarrow \mathfrak{B}(E)$, given by $e \mapsto \mathcal{K}_{e}$, 
is a homomorphism of bounded semilattices which is injective if and only if $E$ is separative.

By an {\em embedding} of a bounded semilattice into a Boolean algebra we mean one that preserves top and bottom elements.

\begin{theorem} Let $E$ be a bounded semilattice.
Then $E$ can be embedded in a Boolean algebra $B$ in such a way that $(E,\wedge)$ is a subsemilattice of $(B,\wedge)$
and each element of $B$ is a join of a finite subset of $E$ if and only if $E$ is separative and compactable.
\end{theorem}
\begin{proof} Let $E$ be a compactable semilattice.
We have proved that $E$ can be embedded in the Boolean algebra $\mathfrak{B}(E)$. 
An element of $\mathfrak{B}(E)$ is a clopen subset.
From the definition of the topology on $\mathsf{U}(E)$ and the fact that it is compact,
such an element can be written as a finite union of clopen sets of the form $\mathcal{K}_{e}$.

Suppose now that $E$ can be embedded in a Boolean algebra $B$ as advertised.
We prove first that $E$ is $0$-disjunctive.
Let $0 \neq f < e$ in $E$.
Then in $B$ there exists $e'$ such that $e' \leq e$ and $f \wedge e' = 0$.
By assumption, we may write $e'$ as a finite join of non-zero elements of $E$.
Let $e_{1}$ be any one of these elements.
Then $e_{1} \wedge f \leq e' \wedge f = 0$ and $e_{1} \leq e$.
It follows that $E$ is $0$-disjunctive.

It remains to show that $E$ is compactable.
Let $e \in E$.
Denote by $\mathcal{K}^{B}_{e}$ the ultrafilters in $B$ containing $e$ 
and by $\mathcal{K}^{E}_{e}$ the ultrafilters in $E$ containing $e$.

We prove that $\mathsf{U}(E)$ is homeomorphic to $\mathsf{U}(B)$.
Let $F$ be an ultrafilter in $E$.
Define 
$$F^{B} = \{b \in B \colon f \leq b \text{ some } f \in F \}.$$
We prove that $F^{B}$ is an ultrafilter in $B$.
Clearly it is a filter.
Observe that $b \in F^{B}$ iff we can write $b = \bigvee_{i=1}^{n} b_{i}$ where $b_{i} \in F$ for some $i$.
To prove this,
suppose therefore that $b \in F^{B}$ but that $b_{i} \notin F$ for all $i$.
By definition $f \leq b$ for some $f \in F$.
Then for each $i$ there exists $f_{i} \in F$ such that $b_{i} \wedge f_{i} = 0$.
Put $f' = \bigwedge_{i} f_{i}$.
Then $f' \in F$ and $b_{i} \wedge f' = 0$ for all $i$.
Thus $f' \wedge b = 0$ and so $f' \wedge f = 0$ which is a contradiction.
The proof of the other direction is clear.

Let $b \in B$ be arbitrary and non-zero and let $b' \in B$ be its complement.
We prove that one of $b$ or $b'$ belongs to $F^{B}$.
Let $b = \bigvee a_{i}$ and $b' = \bigvee c_{j}$ where $a_{i},c_{j} \in E$.
We may assume that $a_{i},c_{j} \notin F$ for all $i$ and $j$.
We may therefore find $f \in F$ such that $f \wedge a_{i} = 0 = f \wedge c_{j}$.
Thus $f \wedge b = 0$ and $f \wedge b' = 0$.
But $1 = b \vee b'$ and $f \wedge 1 \neq 0$,
which gives us a contradiction.
Therefore at least one of the $a_{i} \in F$ or at least one of the $c_{j} \in F$.
It follows that $b \in F^{B}$ or $b' \in F^{B}$.
Thus $F^{B}$ is an ultrafilter in the Boolean algebra $B$.

We have therefore defined a function from $\mathsf{U}(E)$ to $\mathsf{U}(B)$ by $F \mapsto F^{B}$.
Observe that $F^{B} \cap E = F$.
Thus the above function is injective.

Let $G$ be an ultrafilter in $B$.
Put $G^{E} = G \cap E$.
This is non-empty because if $b \in G$ and $b = \bigvee_{i} b_{i}$ where $b_{i} \in E$
then from the properties of ultrafilters in Boolean algebras at least one $b_{i} \in G$.
It is clearly a filter so it remains to show that $G^{E}$ is an ultrafilter in $E$.
Let $e \in E$ be a non-zero element whose meet with every element of $G^{E}$ is non-zero.
Since $G$ is an ultrafilter in $B$ we have that either $e \in G$ or $e' \in G$,
the complement in $B$ of $e$.
Suppose that $e' \in G$.
Then there exists $0 \neq e'' \in G^{E}$ such that $e'' \leq e'$.
But, by assumption, $e \wedge e'' \neq 0$, which is a contradiction.
Thus $e \in G^{E}$ and so $G^{E}$ is an ultrafilter in $E$.

Now $G \cap E \subseteq G$ and so $(G^{E})^{B} \subseteq G$.
But this is an inclusion of ultrafilters and so they must be equal.
We have therefore shown that $(G^{E})^{B} = G$ and so our map is a bijection.

It remains to show that our bijection defines a homeomorphism which will conclude the proof.
Clearly, $\mathcal{K}_{e}^{E}$ is mapped to $\mathcal{K}_{e}^{B}$ and so our map is an open mapping.
It remains to calculate the inverse image under our map of a set of the form $\mathcal{K}_{b}^{B}$.
But
$$\mathcal{K}_{b}^{B}
=
\bigcup_{e \leq b, e \in E} \mathcal{K}_{e}^{B}$$
and the result is now clear.\end{proof}

\begin{theorem} Let $E$ be a separative compactable semilattice
and let $\alpha \colon E \rightarrow B$ be a homomorphism to a Boolean algebra
with the property that the inverse images of ultrafilters are ultrafilters.
Then  there is a unique homomorphism of Boolean algebras $\beta \colon E \rightarrow \mathfrak{B}(E)$
such that $\beta \kappa = \alpha$. 
\end{theorem}
\begin{proof} 
Observe that uniqueness follows from the fact that each element $\mathfrak{B}(E)$
is a join of a finite number of elements from $\kappa (E)$
and that Boolean algebra homomorphisms preserve finite joins.
It thus remains to prove existence.
There is a continuous function $\alpha^{-1} \colon \mathsf{U}(B) \rightarrow \mathsf{U}(E)$ of Boolean spaces
and so a homomorphism of Boolean algebras
$\beta' \colon \mathfrak{B}(E) = \mathsf{B}\mathsf{U}(E) \rightarrow \mathsf{B}\mathsf{U}(B)$.
But $\mathsf{B}\mathsf{U}(B)$ is naturally isomorphic to $B$
and so we have constructed a homomorphism $\beta \colon \mathfrak{B}(E) \rightarrow B$ of Boolean algebras.
In fact, $\beta \kappa = \alpha$ because $\beta (\mathcal{K}_{e}) = \mathcal{K}_{\alpha (e)}$.\end{proof}

Although tangential to the main goals of this paper, we can use the methods of this paper to obtain a characterization of $0$-disjunctive semilattices.
We say that the semilattice $E$ is {\em densely} embedded in the Boolean algebra $B$
if each non-zero element of $B$ lies above a non-zero element of $E$. 

\begin{theorem} A bounded meet semilattice $E$ is $0$-disjunctive if and only if it can be densely embedded in a Boolean algebra $B$.
\end{theorem}
\begin{proof} Let $E$ be a bounded meet semilattice.
From the topological space $\mathsf{U}(E)$ we may still construct a Boolean algebra: namely, the Boolean algebra of all clopen subsets of $\mathsf{U}(E)$.
Denote this Boolean algebra by $B$.
There is an embedding $\kappa \colon E \rightarrow B$ that takes $e$ to $\mathcal{K}_{e}$, which is well-defined because $\mathcal{K}_{e}$ is a clopen set.
Let $X \in B$, an arbitrary clopen set.
Then if $X$ is non-empty there exists a non-zero $e \in E$ such that $\mathcal{K}_{e} \subseteq X$ from the definition of the topology.

Conversely, suppose that $E$ is a semilattice embedded in the way advertised in the Boolean algebra $B$.
Let $f < e$ in $E$.
Then since $B$ is a Boolean algebra there exists a non-zero element $f' \leq e$ such that $f \wedge f' = 0$.
But, by assumption, there exists $0 \neq f'' \leq f'$ where $f''\in E$.
Clearly $f'' \leq e$ and $f'' \wedge f = 0$.
\end{proof}

To show that the theory is not vacuous, we shall construct some examples of compactable bounded meet semilattices which are not Boolean algebras. 
The examples we give will not be surprising from the point of view of the theory of $C^{\ast}$-algebras \cite{Lenz,P};
it is rather the context which is novel.

Let $E$ be a bounded semilattice. 
Given $e,f \in E$ we say that $e$ {\em covers} $f$ if $e > f$ and there is no $g \in E$ such that $e > g > f$.
For each $e \in E$ define $\hat{e}$ to be the set of elements of $E$ that are covered by $e$.
A semilattice is said to be {\em pseudofinite} if whenever $e > f$ there exists $g \in \hat{e}$
such that $e > g > f$, and for which the sets $\hat{e}$ are always finite.
This definition was used by Munn \cite{M}.

A semilattice is said to be {\em unambiguous} if for all non-zero elements $e$ and $f$
we have that $e \wedge f \neq 0$ implies that either $e \leq f$ or $f \leq e$.
In other words, the posets $e^{\uparrow}$, where $e \neq 0$, are linearly ordered.
This definition is a special case of one that can be made for any semigroup with zero
although the standard terminology  is {\em unambiguous except at zero}.
Such semigroups play an important role in those parts of semigroup theory motivated by the Krohn-Rhodes theorem \cite{B1,B2}.
Unambiguous {\em inverse} monoids were first discussed in my paper \cite{L3}.

\begin{proposition} Let $E$ be a bounded $0$-disjunctive semilattice which is 
pseudofinite, unambiguous
and in which $e^{\uparrow}$ is finite for each non-zero element $e$.
Then $E$ is compactable.
\end{proposition}
\begin{proof} Define a function $\lambda \colon E \rightarrow \mathbb{N} \cup \{\infty \}$ 
by $\lambda (e) = \left| e^{\uparrow} \right|$ if $e \neq 0$ and $\lambda (0) = \infty$.
We call $\lambda (e)$ the {\em level} of $e$.
This function has some useful properties:
\begin{itemize}

\item Observe that $f \leq e$ iff $e^{\uparrow} \subseteq f^{\uparrow}$.
Thus  $f < e$ iff $e^{\uparrow} \subset f^{\uparrow}$.
It follows that $f < e$ implies that $\lambda (f) > \lambda (e)$.

\item If $\lambda (e) = \lambda (f)$ then either $e = f$ or $e$ and $f$ are orthogonal.
To see why, suppose that $e \wedge f \neq 0$.
Then either $e < f$ or $f < e$ or $e = f$.
Suppose that $e \neq f$. 
Then either $\lambda (e) > \lambda (f)$ or $\lambda (f) > \lambda (e)$, which is a contradiction.
\end{itemize}
We show that the trapping condition holds by using induction.
Let $e \in E$ and let $f < e$ such that $\lambda (f) = \lambda (e) + 1$.
By pseudofiniteness, the order ideal $e^{\downarrow}\setminus \{ e\}$ is generated by a finite set $\{e_{1}, \ldots, e_{m}\}$.
By unambiguity, we may assume that the elements $e_{i}$ are pairwise orthogonal and all have level one more than that of $e$.
By pseudofiniteness, 
from $f < e$ we must have $f = e_{i}$ for some $i$. Without loss of generality we may assume that $f = e_{1} $.
Observe that by pseudofinitness $e \rightarrow (e_{2}, \ldots, e_{m},f)$ and the trapping condition holds. 
We have therefore proved that the trapping condition holds for all pairs $f < e$ which differ by one level.

Suppose now that the trapping condition holds for all $f < e$ such that $\lambda (f) = \lambda (e) + r$.
We prove that it holds for all $f < e$ such that $\lambda (f) = \lambda (e) + r+1$.
There exists an $f' \in E$ such that $f < f' < e$ where $\lambda (f) = \lambda (f') + 1$.
Thus by our induction hypothesis, the trapping condition holds for $f' < e$ and so there exist
elements $e_{i} \leq e$ and orthogonal to $f'$ such that $e \rightarrow (e_{1}, \ldots, e_{m},f')$.

We have that $f < f'$ and these elements differ in one level so we may find elements $f_{1}, \ldots, f_{n} \leq f'$ orthogonal to $f$ such that
$f' \rightarrow (f_{1}, \ldots, f_{n},f)$.
We claim that
$$e \rightarrow (e_{1}, \ldots, e_{m},f_{1}, \ldots, f_{n},f).$$ 
Let $0 \neq x \leq e$.
Then either $x \wedge e_{i} \neq 0$  for some $i$ or $x \wedge f' \neq 0$.
If the latter then $0 \neq x \wedge f' \leq f'$.
Thus either $x \wedge f' \wedge f_{j} \neq 0$ for some $j$ or $x \wedge f' \wedge f \neq 0$.
It follows that  
$x \wedge e_{i} \neq 0$  for some $i$
or
$x \wedge f_{j} \neq 0$ for some $j$
or
$x \wedge f \neq 0$.
It follows that $e \rightarrow (e_{1}, \ldots, e_{m},f_{1}, \ldots, f_{n},f)$.
Now $e_{i} \in (f')^{\ast}$ and $f_{j} \in f^{\ast}$ and $f < f'$.
Hence $e_{i},f_{j} \in f^{\ast}$ and so the trapping condition holds.\footnote{This part of the proof uses the fact
that the relation $\rightarrow$ is `transitive'. See Proposition~5.6 of \cite{Lenz}.}
\end{proof}

Concrete examples of semilattices satisfying the conditions of the above proposition may easily be constructed.
Let $G$ be a directed graph with a distinguished vertex $t$
such that for each vertex $v$ of $G$ there is a path in $G$ that starts at $v$ and ends at $t$.
We shall say that such a graph is {\em rooted}.
We denote the free category on $G$ by $G^{\ast}$.
We are interested in the set of paths that start at $t$.
Let $E$ be the set of such paths together with a zero element.
Define $e \leq f$ iff $e = fg$ for some element $g \in G^{\ast}$;
in other words, the prefix ordering.
It is easy to check that with respect to this order, 
$E$ becomes a bounded semilattice which is unambiguous and in which the set $e^{\uparrow}$ is finite
for all non-zero elements $e$.

\begin{lemma} With the above definitions, we have the following.
\begin{enumerate}

\item The semilattice $E$ is $0$-disjunctive iff the in-degree of each vertex is either zero or at least 2.

\item  The semilattice $E$ is pseudofinite if and only if the in-degree of each vertex is finite.

\end{enumerate}
\end{lemma}
\begin{proof} 
(1) Suppose that $E$ is $0$-disjunctive.
Let $v$ be any vertex and let $x$ be a path from $v$ to $t$.
Suppose that the in-degree of $v$ is not zero.
Then there is at least one edge $w$ into $v$.
It follows that $xw \leq x$.
By assumption, there exists $z \leq x$ such that $z$ and $xw$ are orthogonal.
Now $z = xp$ for some non-empty path $p$.
It follows that $w$ is not a prefix of $p$ and so there is at least one other edge coming into the vertex $v$.
 
Suppose now that the in-degree of each vertex is either zero or at least two.
Let $y < x$ where $y  = xp$.
Let $x$ start from the vertex $v$.
Since $p$ is a non-empty path that starts at $v$ it follows that there is at least one other edge $w$ coming into $v$
that differs from the first edge of $p$.
Thus $xw \leq x$ and $xw$ and $y$ are orthogonal.

(2) This is clear.

\end{proof}

\begin{example}
{\em 
The simplest example of a compactable semilattice of the above type which is not a Boolean algebra 
occurs in the case where the directed graph consists of one vertex and two loops:
in other words, the semilattice $E_{2}$ obtained from the free monoid $A^{\ast}$ on two generators $A = \{a,b \}$ 
via the prefix ordering and the adjoining of a zero.
We put 1 equal to the empty string.
Observe that $E_{2}$ is not a Boolean algebra.
To see why, observe that we have $aa \leq 1$ and that $b,ab$ are both orthogonal to $aa$.
But the only element above both $b$ and $ab$ is 1.
Thus $aa$ has no Boolean complement.

We now calculate $\mathfrak{B}(E_{2})$.
Denote by $A^{\omega}$ the set of all right-infinite strings over the alphabet $A$.
The ultrafilters in $E_{2}$ can be identified with the elements of $A^{\omega}$ as follows.
Let $w \in A^{\omega}$ and put 
$$F_{w} = \{x \in A^{\ast} \colon x \text{ is a prefix of } w \}.$$
Then this is an ultrafilter in $E_{2}$ and every ultrafilter in $E_{2}$ is of this form.
We may therefore identify $\mathsf{U}(E_{2})$ with $A^{\omega}$.
For each $x \in A^{\ast}$, we may identify the set $\mathcal{K}_{x}$ with all those infinite strings that begin with $x$.
This is a base for the topology we have defined on $\mathsf{U}(E_{2})$ but it is also the base for the product topology on $A^{\omega}$.
It follows that  $\mathfrak{B}(E_{2})$ is isomorphic to the Boolean algebra of clopen subsets of $A^{\omega}$.
Such sets are of the form $XA^{\omega}$ where $X \subseteq A^{\ast}$ is a finite subset \cite{PP}.
}
\end{example}

Finally, there is an evident connection between this paper and those by Hughes \cite{H1,H2}.
The semilattices we constructed above are examples of $\mathbb{R}$-trees since they are rooted locally finite simplicial trees.
For such $\mathbb{R}$-trees there is then a dictionary between Hughes' geometry and our algebra:
{\em geodesically complete} implies that there are no $0$-minimal idempotents;
{\em ends} correspond to ultrafilters;
{\em cut sets} correspond to maximal prefix codes.
There are also parallels between \cite{H2} and our paper \cite{Lawson}:
this is not surprising given the well-known connections between \'{e}tale topological groupoids and inverse semigroups \cite{P}.
However, it suggests an interesting avenue for future research.


\end{document}